\documentclass{amsart}

\usepackage{amsmath, graphicx}
\usepackage{amsfonts}
\usepackage{amssymb}
\usepackage{hyperref}

\newcommand{\summ}{\mathop{{\sum}^{\star}}}

\numberwithin{equation}{section}

\newcommand{\sums}{\mathop{{\sum}^{*}}}

\newcommand{\sym}{\mathrm{sym}^2}

\newtheorem{theorem}{Theorem}[section]

\newtheorem{lemma}[theorem]{Lemma}

\begin{document}

\title{$L^4$-norms of Hecke newforms of large level}
\author{Jack Buttcane and Rizwanur Khan}
\address{Mathematisches Institut, Georg-August Universit\"{a}t G\"{o}ttingen, Bunsenstrasse 3-5, D-37073 G\"{o}ttingen, Germany}
\email{buttcane@uni-math.gwdg.de, rrkhan@uni-math.gwdg.de}
\thanks{}

\begin{abstract} We prove a new upper bound for the $L^4$-norm of a holomorphic Hecke newform of large fixed weight and prime level $q\to \infty$. This is achieved by proving a sharp mean value estimate for a related $L$-function on $GL(6)$.
\end{abstract}

\maketitle

\section{Introduction}

Modular forms, particularly cusp forms, are very special functions on the upper half complex plane which arise almost everywhere in mathematics. A natural way to understand a cusp form is to study its $L^p$-norms, for in principle a function can be recovered from the knowledge of its moments. In recent years, this has been an exciting topic for which even a 2010 Fields Medal was awarded, to Lindenstrauss. Settling an important case of the Quantum Unique Ergodicity conjecture, Lindenstrauss \cite{lin} and Soundararajan \cite{sou} showed that the $L^2$-norm of a Hecke Maass cusp form of large laplacian eigenvalue restricted to a finite region of the complex upper half plane depends only on the area of the region. In other words, the $L^2$-mass is equidistributed. The analogous problem for holomorphic Hecke cusp forms of large weight, the Rudnick-Sarnak conjecture, was solved by Holowinsky and Soundararajan \cite{holsou}. Nelson \cite{nel} proved that a holomorphic Hecke newform of large square-free level has equidistributed $L^2$-mass. This version of the QUE conjecture in the level aspect was posed by Kowalski, Michel, and VanderKam.

What can be said about {\it higher} $L^p$-norms in the level aspect? Let $\mathcal{B}_k^{\text{\it{new}}}(q)$ denote the set of $L^2$-normalized holomorphic Hecke newforms of level $q$, weight $k$, and trivial nebentypus, and let $f\in \mathcal{B}_k^{\text{\it{new}}}(q)$.  Blomer and Holowinsky \cite{blohol} were the first to establish a non-trivial upper bound for the $L^\infty$-norm. This was later improved by Harcos and Templier\footnote[1]{These authors actually proved their results for Hecke Maass forms, but very similar arguments would work for holomorphic Hecke newforms as well.} \cite{hartem}, who showed that
\begin{equation}
\| f \|_\infty \ll q^{-1/6+\epsilon} 
\end{equation}
for any $\epsilon>0$, where the implied constant depends on $\epsilon$ and $k$. (We adopt this convention throughout the paper, and also that $\epsilon$ may not be the same arbitrarily small positive constant from one occurrence to the next.) As a consequence we have that
\begin{equation}
\label{indirect} \| f \|_4^4 \ll  q^{-1/3 +\epsilon},
\end{equation}
which is an improvement over the ``trivial'' bound
\begin{equation}
\| f \|_4 \ll q^\epsilon .
\end{equation}
A more direct way to get a handle on the $L^4$-norm is through $L$-functions. Blomer \cite[sec. 1]{blo2} commented that ``it seems to be very hard to improve'' (\ref{indirect}) in this way, but nevertheless Liu, Masri, and Young\footnotemark[1] 
\cite{liumasyou} succeeded in the case that $q$ is prime and achieved the bound
\begin{equation}
\| f \|_4^4 \ll q^{-1/2 +\epsilon}.
\end{equation}
In this paper, we give a further improvement. Let $q$ denote a prime throughout the rest of the paper.
\begin{theorem} \label{main}
Let $\epsilon>0$. There exists $k_\epsilon>0$ such that for even $k>k_{\epsilon}$ the following result holds. Let $q$ be prime and $f\in \mathcal{B}_{k}^{\text{\it{new}}}(q)$. Suppose that the subconvexity bound
\begin{equation}
\label{subconvex} L(\tfrac{1}{2}, g) \ll_{k,\epsilon} q^{1/4-\delta+\epsilon}
\end{equation}
holds for all $g\in\mathcal{B}_{2k}^{\text{\it{new}}}(q)$, where $L(s,g)$ is defined in (\ref{assoc}) and $\delta>0$. Then we have that
\begin{align}
\| f \|_4^4 \ll_{k,\epsilon} q^{-3/4-\delta+\epsilon}.
\end{align}
\end{theorem}
\noindent We have concentrated only on arriving at a new bound in the level aspect, and have not tried to optimize $k_\epsilon$. The best known subconvexity bound, due to Duke, Friedlander, and Iwaniec \cite{dfi}, is given by $\delta=1/192$. It is conjectured of course that we can take $\delta=1/4$; this would yield the best expected bound for $\| f \|_4$.

The proof of Theorem \ref{main} rests on a new mean value estimate for an $L$-function on $GL(6)$, which is of interest in its own right. 
\begin{theorem}\label{meanval} Let $\epsilon>0$. There exists $k_\epsilon>0$ such that for even $k>k_\epsilon$, $q$ prime, and $f\in\mathcal{B}_k^{\text{\it{new}}}(q)$, we have that
\begin{equation}
\label{meanvalue} \sum_{g \in\mathcal{B}_{2k}^{\text{\it{new}}}(q)} \frac{L(\tfrac{1}{2}, \sym f \times g)}{L(1,\sym g)} \ll_{k,\epsilon} q^{1+\epsilon}. 
\end{equation}
\end{theorem}
\noindent After noting (\ref{lower}), we see that the estimate above is consistent with the Lindel\"{o}f hypothesis. Further, using Lapid's \cite{lap} result that 
\begin{equation}
\label{lapid} L(\tfrac{1}{2}, \sym f \times g)\ge 0,
\end{equation}
dropping all but one term recovers the convexity bound for these central values. 

The novelty of our paper which facilitates the proof of Theorem \ref{meanval} is a type of $GL(3)$ Voronoi summation formula for non-trivial level. The $GL(3)$ Voronoi summation formula for level one, due to Miller and Schmid \cite{schmil}, is a valuable tool in analytic number theory which has been used extensively in the subject. Although Ichino and Templier \cite{ichtem} have established a more general formula, it does not cover the case of interest to us. We will need to analyze a sum of the type
\begin{equation}
\label{type} \sum_{n<N} A_f(n,1)e\Big( \frac{nh}{cq} \Big)
\end{equation}
where  $A_f(n,1)$ are the coefficients of $L(s,\sym f)$ as given in (\ref{symsq}), $c$ is an integer less than $q^{\epsilon}$, and $h$ is an integer coprime to $cq$. It would be very desirable to obtain a practical formula for the above sum for general values of $c$, but this seems to be a challenging task. Our idea to get around this problem in the current situation is to realize that it is enough to develop a formula for the telling case $c=1$ but to not fully work out the contribution coming from $c>1$. Since $c<q^\epsilon$, this approach is good enough. 

\begin{lemma}\label{voronoi} Let $A_f(n,1)$ be the coefficients of $L(s,\sym f)$ as given in (\ref{symsq}). Let $\phi(x)$ be a smooth function compactly supported on the positive real numbers.
 For $(h,q)=1$, we have that
\begin{multline}
\label{vor} \sum_{\substack{n\ge 1\\ (n,q)=1}} A_f(n,1) e\Big( \frac{nh}{q} \Big) \phi(n) =  \frac{q}{2} \sum_{\alpha=\pm 1}  i^{\frac{1+\alpha}{2}} \sum_{\substack{n\ge 1\\(n,q)=1}} \frac{A_f(n,1)}{n}\Big(S(-n\overline{h},1;q)+\alpha S(n\overline{h},1;q)\Big)\Phi_\alpha \Big(\frac{n}{q^3}\Big)\\ -\sum_{\substack{n\ge 1 \\ (n,q)=1}} \frac{A_f(n,1)}{q-1} \Big(\phi(n) + \frac{iq}{n} \Phi_1 \Big(\frac{n}{q^3}\Big)\Big),
\end{multline}
where $\overline{h}$ is the multiplicative inverse of $h$ modulo $q$ and for $\widetilde{\phi}(s)$ the Mellin transform of $\phi$ and $H_{\alpha}(s)$ as in (\ref{g2}), we define
\begin{equation}
\Phi_\alpha(x) =  \frac{1}{2\pi i}  \int_{(\sigma)} x^{-s}  \frac{H_{\alpha}(1-s)}{H_{\alpha}(s)} \widetilde{\phi}(s) \ ds
 \end{equation}
for any $x,\sigma>0$.
\end{lemma}

It would be interesting to prove an analogue of Theorem \ref{meanval} for Hecke Maass cusp forms of level $q$, and we expect that very similar arguments would work. This would lead to a new bound for the $L^4$-norm of Maass forms and an improvement of the subconvexity bound \cite[Corollary 1.2]{liumasyou}. That a result like Theorem \ref{meanval} would improve their subconvexity bound was noted by Liu, Masri, and Young in \cite[section 6]{liumasyou}.

\subsection{Sketch} We now give a rough outline of the proof of Theorem \ref{meanval}. 

By the method of approximate functional equations, we write
\begin{equation}
L(\tfrac{1}{2}, \sym f \times g) \approx \sum_{n< q^{2+\epsilon}} \frac{A_f(n,1)\lambda_g(n)}{n^{1/2}},
\end{equation}
 where $\lambda_g(n)$ are the coefficients of $L(s,g)$ as given in (\ref{assoc}). Thus,
\begin{equation}
\sum_{g \in\mathcal{B}_{2k}^{\text{\it{new}}}(q)} \frac{L(\tfrac{1}{2}, \sym f \times g)}{L(1,\sym g)} \approx
\sum_{n <q^{2+\epsilon}} \frac{A_f(n,1)}{n^{1/2}}  \sum_{g \in\mathcal{B}_{2k}^{\text{\it{new}}}(q)}\frac{\lambda_g(n)}{L(1,\sym g)}.
 \end{equation}
 By the Petersson trace formula, this is roughly
 \begin{equation}
q\left( 1+ \sum_{n <q^{2+\epsilon}} \frac{A_f(n,1)}{n^{1/2}} \sum_{c\ge 1} J_{2k-1}\Big(\frac{4\pi\sqrt{n}}{cq} \Big)\frac{S(n,1;cq)}{cq}\right).
 \end{equation}
Now for large enough $k$, the $J$-Bessel function $J_{2k-1}(x)$ decays very quickly as $x\to 0$. Therefore it suffices to consider only $n>q^{2-\epsilon}$ and $c<q^{\epsilon}$. To fix ideas, we restrict ourselves in this sketch to the case $q^2<n<2q^2$ and $c=1$. In these ranges, $n^{1/2}\approx q$ and the $J$-Bessel function is roughly constant. We therefore have to show that
\begin{equation}
\sum_{q^2<n<2q^2} A_f(n,1)S(n,1;q) \ll q^{2+\epsilon}.
\end{equation}

Opening the Kloosterman sum, we have to show that
\begin{equation}
\label{sketch} \sums_{h \bmod q} e\Big(\frac{\overline{h}}{q}\Big) \sum_{q^2<n<2q^2} A_f(n,1)e\Big(\frac{nh}{q}\Big) \ll q^{2+\epsilon},
\end{equation}
where $*$ restricts the sum to primitive residue classes. We develop a Voronoi summation formula with the following shape:
\begin{equation}
\sum_{q^2<n<2q^2} A_f(n,1)e\Big(\frac{nh}{q}\Big) \approx q\sum_{m< q^{1+\epsilon}} \frac{A_f(m,1)}{m} S(m,\overline{h};q).
\end{equation}
Inserting this into (\ref{sketch}), we have to show that
\begin{equation}
\sum_{m< q^{1+\epsilon}} \frac{A_f(m,1)}{m} \sums_{h \bmod q} e\Big(\frac{\overline{h}}{q}\Big)S(m,\overline{h};q) \ll q^{1+\epsilon}.
\end{equation}
Evaluating the $h$-sum, we must show that
\begin{equation}
\sum_{m< q^{1+\epsilon}} \frac{A_f(m,1)}{m} e\Big(\frac{-m}{q}\Big) \ll q^{\epsilon}.
\end{equation}
The required bound now follows by estimating trivially. Proving any further cancelation in the sum above seems like a very difficult problem, as the length of the sum is the square root of the conductor of $L(s,\sym f)$.

To obtain the Voronoi formula, we write the exponential in terms of Dirichlet characters to get that
\begin{equation}
\sum_{n\ge 1} A_f(n,1)e\Big(\frac{nh}{q}\Big)\phi\Big(\frac{n}{q^2}\Big) \approx \frac{1}{q} \sum_{\substack{\chi \bmod q\\ \chi\neq 1}} \tau(\chi)  \sum_{n\ge 1} A_f(n,1) \overline{\chi}(nh)\phi\Big(\frac{n}{q^2}\Big),
\end{equation}
where $\phi(x)$ is a smooth bump function supported on $1<x<2$ and $\tau(\chi)$ is the Gauss sum. In terms of $L$-functions, this equals
\begin{equation}
\frac{1}{q}\sum_{\substack{\chi \bmod q\\ \chi\neq 1}} \tau(\chi) \overline{\chi}(h) \frac{1}{2 \pi i}  \int_{2-\infty}^{2+i\infty} L(s, \sym f \times \overline{\chi}) q^{2s}    \widetilde{\phi}(s) ds,
\end{equation}
where $\widetilde{\phi}$ is the Mellin transform of $\phi$. The summation formula comes from an application of the functional equation of $L(s, \sym f \times \overline{\chi})$, which is available from the work of Li \cite{li}.

\section{Background}

\subsection{Cusp forms}

Let $\mathbb{H}$ denote the upper half complex plane and $\Gamma_0(q)$ the usual congruence subgroup of $SL_2(\mathbb{Z})$. Let $S_k(q)$ denote the space of cusp forms of weight $k$ and trivial nebentypus for $\Gamma_0(q)$. This space is equipped with the inner product 
\begin{equation}
\langle f_1,f_2 \rangle =  \int_{\Gamma_0(q)\backslash \mathbb{H}} y^{k/2} f_1(z) \overline{y^{k/2} f_2(z)} \frac{dx dy}{y^2}
\end{equation}
and the $L^p$-norms
\begin{align}
&\| f \|_p = \Big( \int_{\Gamma_0(q)\backslash \mathbb{H}} |y^{k/2}f(z)|^p \frac{dx dy}{y^2} \Big)^{1/p}, \\
&\| f \|_\infty = \sup \{  |y^{k/2}f(z)| : z\in \Gamma_0(q)\backslash \mathbb{H} \}.
\end{align}

By \cite[Proposition 2.6]{ils} and \cite[sec. 2]{blo2}, there exists an orthonormal basis $\mathcal{B}_{k}^\text{\it{new}}(q)\cup \mathcal{B}_{k}^\text{\it{old}}(q)$ for $S_k(q)$, where every $f\in \mathcal{B}_{k}^\text{\it{old}}(q)$ is an  oldform with
\begin{equation}
\label{oldform} \| f \|_\infty \ll q^{-1/2}.
\end{equation}
We have that
\begin{equation}
\label{dim} \dim S_k(q) \sim |\mathcal{B}_{k}^\text{\it new}(q)| \sim \frac{q(k-1)}{12}
\end{equation}
as $q\to \infty$.

\subsection{$L$-functions} Every $f\in \mathcal{B}_{k}^{\text{\it new}}(q)$ is an eigenfunction of the Hecke operators $T_n$. Say
\begin{equation}
T_n f = n^{\frac{k-1}{2}} \lambda_f(n) f
\end{equation}
for some real numbers $ \lambda_f(n)$ satisfying 
\begin{equation}
\label{fricke} \lambda_f(q)\ll q^{-1/2},
\end{equation}
Deligne's bound $\lambda_f(n)\ll n^\epsilon$, and the multiplicative relation
\begin{equation}
\label{heckemult} \lambda_f(n)\lambda(m) = \sum_{\substack{d|(n,m) \\ (d,q)=1}} \lambda_f\Big(\frac{nm}{d^2}\Big).
\end{equation}

The $L$-function associated to $f$ equals
\begin{equation}
\label{assoc} L(s,f) = \sum_{n\ge 1} \frac{\lambda_f(n)}{n^s}
\end{equation}
for $\Re(s)>1$ with analytical continuation to the rest of the complex plane. By the work of Guo \cite{guo}, we know that the central value is non-negative:
\begin{equation}
\label{guo} L(\tfrac{1}{2},f) \ge 0.
\end{equation}

We will also need to work with $GL(1)$ and $GL(2)$ twists of the symmetric-square $L$-function, given for $\Re(s)>1$ by
\begin{equation}
\label{symsq} L(s,\sym f) = \Big(1-\frac{1}{q^{s+1}}\Big)^{-1}  \sum_{\substack{n\ge 1\\ (n,q)=1}} \frac{A_f(n,1)}{n^s},
\end{equation}
where
\begin{equation}
\label{gl3coeff} A_f(n,1)= \sum_{d_1^2 d_2| n} \lambda_f(d_2^2)
\end{equation}
for $(n,q)=1$. At the edge of the critical strip have the bounds
\begin{align}
\label{lower}  C_\epsilon q^{-\epsilon} < L(1,\sym f) \ll q^{\epsilon},
\end{align}
where $C_\epsilon$ is some positive constant depending on $\epsilon$. The lower bound is due to Goldfeld, Hoffstein, and Lieman \cite{golhoflie}. The constant $C_\epsilon$ is ineffective if $f$ is dihedral, but in this case even better bounds for the $L^4$-norm may be available by arguments such as in \cite{luo}.

We will consider for $(c,q)=1$ and $\chi$ a primitive Dirichlet character of modulus $cq$, the twist
\begin{align}
\label{twist} L(s,\sym f \times \chi) = \sum_{n\ge 1} \frac{A_f(n,1)\chi(n)}{n^s}
\end{align}
for $\Re(s)>1$.
Its functional equation is given in section \ref{sect:fetwist}.

Let $g \in B_{2k}^{\text{\it new}}(q)$ and
\begin{equation}
\label{mult} A_f(n,m) = \sum_{d|(n,m)}\mu(d) A_f\Big( \frac{n}{d} , 1 \Big) A_f\Big( \frac{m}{d} , 1 \Big)
\end{equation}
for $(nm,q)=1$. By \cite[section 3.1]{wat}, the $GL(2)$ twist
\begin{equation}
 \label{series} L(s,\sym f \times g) = \Big(1-\frac{\lambda_g(q)}{q^{s}}\Big)^{-1} \Big(1-\frac{\lambda_g(q)}{q^{s+1}}\Big)^{-1} \sum_{\substack{n,m\ge 1\\(nm,q)=1}} \frac{A_f(n,m) \lambda_g(n) }{ n^s m^{2s} }
 \end{equation}
for $\Re(s)>1$ continues to an entire function and has the functional equation
\begin{equation}
\label{fe} q^{2s} G(s)  L(s,\sym f \times g)  = q^{2(1-s)} G(1-s) L(1-s,\sym f \times g),
\end{equation}
where
\begin{equation}
G(s) =  \pi^{-3s} \Gamma(s+2k-\tfrac{3}{2}) \Gamma(s+k-\tfrac{1}{2}) \Gamma(s+\tfrac{1}{2}) .
\end{equation}
We will use an approximate functional equation to get a handle on the central values.
\begin{lemma}We have
\begin{align}
\label{afe} L(\tfrac{1}{2}, \sym f \times g) = 2\sum_{\substack{n, m \ge 1 \\ (nm,q)=1}}  \frac{A_f(n,m) \lambda_g(n) }{ n^{1/2} m } \sum_{r_1,r_2\ge 0}\Big( \frac{\lambda_g(q)}{q^{1/2}} \Big)^{r_1} \Big( \frac{\lambda_g(q)}{q^{3/2}} \Big)^{r_2} V\Big( \frac{ nm^2}{ q^{2-r_1-r_2}}\Big),
\end{align}
where
\begin{equation}
V(x) = \frac{1}{2\pi i}  \int_{(\sigma)} x^{-s}  \frac{G(\tfrac{1}{2} +s)}{G(\tfrac{1}{2})} \frac{ds}{s}
\end{equation}
satisfies
\begin{equation}
\label{v} {V}^{(\ell)}(x) \ll_{\ell,\sigma} x^{-\ell-\sigma}
\end{equation}
for $x,\ell,\sigma>0$. Thus the sum in (\ref{afe}) is essentially supported on $nm^2<q^{2-r_1-r_2+\epsilon}$.
\end{lemma}
\proof
By Cauchy's theorem, we have that
\begin{align}
\label{afeproof} L(\tfrac{1}{2}, \sym f \times g) = &\frac{1}{ 2\pi i} \int_{\big(\tfrac{3}{4}\big)} \frac{L(\tfrac{1}{2}+s,\sym f \times g)q^{2s} G(\tfrac{1}{2}+s)}{G(\tfrac{1}{2})} \frac{ds}{s}\\
\nonumber - &\frac{1}{ 2\pi i} \int_{\big(-\tfrac{3}{4}\big)} \frac{L(\tfrac{1}{2}+s,\sym f \times g)q^{2s} G(\tfrac{1}{2}+s)}{G(\tfrac{1}{2})} \frac{ds}{s}.
\end{align}
Applying the functional equation to the integrand in the second line of (\ref{afeproof}), we get that
\begin{align}
L(\tfrac{1}{2}, \sym f \times g) = \frac{1}{ \pi i} \int_{\big(\tfrac{3}{4}\big)} \frac{L(\tfrac{1}{2}+s,\sym f \times g)q^{2s} G(\tfrac{1}{2}+s)}{G(\tfrac{1}{2})} \frac{ds}{s}.
\end{align}
On this line of integration we may use (\ref{series}) and the Taylor expansion of the Euler factors at $q$. Doing so, we arrive at (\ref{afe}) with $\sigma=3/4$. Now the line of integration may be moved to any $\sigma>0$.
 \endproof

\subsection{Trace formula}
We will need the following trace formula over newforms, implied by \cite[Propositions 2.1 and 2.8]{ils}. For $(n,q)=1$, we have that 
\begin{align}
\label{trace} \frac{12\zeta(2)}{q(k-1)}\sum_{f\in B_k^{\text{\it new}}(q)} \frac{\lambda_f(n)}{L(1,\sym f)} = \delta(n) + 2\pi i^k \sum_{c\ge 1} \frac{S(n,1;cq)}{cq}J_{k-1}\Big( \frac{4\pi\sqrt{n} }{cq} \Big)+O\Big(\frac{n^{\epsilon}}{q} \Big),
\end{align}
where $\delta(n)=1$ if $n=1$ and $\delta(n)=0$ otherwise, $S(n,1;cq)$ is a Kloosterman sum, and $J_{k-1}(x)$ is the $J$-Bessel function, which satisfies
\begin{align}
\label{decay} J_{k-1}(x) \ll \min \{ x^{k-1}, x^{-1/2} \}
\end{align}
for $x>0$.
We have the very basic uniform bound
\begin{equation}
 J_{k-1}^{(\ell)}(x) \le 1
\end{equation}
for all $x>0$ and $\ell \ge 0$.

\section{Reduction of Theorem \ref{main}}

The goal of this section is to relate the $L^4$-norm to $L$-functions and reduce Theorem \ref{main} to Theorem \ref{meanval}. Note that if $f\in B_k^\text{\it new}(q)$ then $f^2 \in S_{2k}(q)$. Thus by Parseval's theorem and (\ref{oldform}-\ref{dim}), we have
\begin{equation}
\| f \|_4^4 = \| f^2 \|_2^2 = \sum_{g\in B_{2k}^\text{\it new}(q)\cup B_{2k}^\text{\it old}(q)} |\langle f^2 ,g \rangle|^2 = \sum_{g\in B_{2k}^\text{\it new}(q)} |\langle f^2 ,g \rangle|^2 + O(q^{-1}).
\end{equation}
By Watson's formula \cite[section 4.1]{wat} (see also \cite[sec. 5]{blo2}) for $|\langle f^2 ,g \rangle|^2$, we have that
\begin{align}
\label{watsonapplied}\| f \|_4^4 \ll \frac{1}{q^{2}} \sum_{g\in  B_{2k}^\text{\it new}(q)} \frac{L(\tfrac{1}{2}, f\times f\times \overline{g})}{L(1,\sym f)^2L(1, \sym g)}+ O(q^{-1}),
\end{align}
where
\begin{align}
L(s, f\times f\times \overline{g}) = L(s, g)L(s,\sym f \times g)
\end{align}
is a triple product $L$-function (the complex conjugation can be dropped because $\lambda_g(n)\in \mathbb{R}$). Now by (\ref{subconvex}, \ref{lapid}, \ref{guo}, \ref{lower}), we have that
\begin{equation}
\| f \|_4^4 \ll q^{-7/4-\delta+\epsilon} \sum_{g\in  B_{2k}^\text{\it new}(q)} \frac{L(\tfrac{1}{2},\sym f \times g)}{L(1,\sym g)} + O(q^{-1}).
\end{equation}
Thus Theorem \ref{main} follows from Theorem \ref{meanval}. The same strategy was used in \cite{blokhayou} to study the $L^4$-norm in terms of the weight $k$. Notice how this differs from the approach in \cite{liumasyou}: there, the Cauchy-Schwarz inequality is applied to (\ref{watsonapplied}) and the problem is reduced to studying the second power moments of $L(\tfrac{1}{2}, g)$ and $L(\tfrac{1}{2},\sym f \times g)$, for Maass forms.

We remark that the Lindel\"{o}f bound for the triple product $L$-function above would imply the best expected bound $\| f \|_4 \ll q^{-1/4+\epsilon}$.

\section{Proof of Theorem \ref{meanval}}

By (\ref{afe}), we have that
\begin{multline}
 \frac{1}{q} \sum_{g \in\mathcal{B}_{2k}^{\text{\it{new}}}(q)} \frac{L(\tfrac{1}{2}, \sym f \times g)}{L(1,\sym g)}\\
 =\frac{2}{q}  \sum_{g \in\mathcal{B}_{2k}^{\text{\it{new}}}(q)}\sum_{\substack{n, m \ge 1 \\ (nm,q)=1}}  \frac{A_f(n,m) } { n^{1/2} m }\frac{\lambda_g(n)} {L(1, \sym g)} \sum_{r_1,r_2\ge 0}  \Big( \frac{\lambda_g(q)}{q^{1/2}} \Big)^{r_1} \Big( \frac{\lambda_g(q)}{q^{3/2}} \Big)^{r_2} V\Big( \frac{ nm^2}{ q^{2-r_1-r_2}}\Big).
 \end{multline}
 Now by (\ref{fricke}), this equals
\begin{equation}
2 \sum_{\substack{n, m \ge 1 \\ (nm,q)=1}}  \frac{A_f(n,m)}{ n^{1/2} m } V\Big( \frac{ nm^2}{ q^2}\Big)  \frac{1}{q}  \sum_{g \in\mathcal{B}_{2k}^{\text{\it{new}}}(q)} \frac{\lambda_g(n)}{L(1, \sym  g)} +O(q^{\epsilon}).
\end{equation}
Thus on applying (\ref{trace}), we have that
\begin{equation}
\label{aftertrace} \frac{1}{q} \sum_{g \in\mathcal{B}_{2k}^{\text{\it{new}}}(q)} \frac{L(\tfrac{1}{2}, \sym f \times g)}{L(1,\sym g)} \ll \sum_{c\ge 1} \sum_{\substack{n, m \ge 1 \\ (nm,q)=1}}  \frac{A_f(n,m) } { n^{1/2} m }\frac{S(n,1;cq)}{cq} J_{2k-1}\Big( \frac{4\pi\sqrt{n} }{cq}  \Big)V\Big( \frac{ nm^2}{ q^2}\Big)+O(q^{\epsilon}).
\end{equation}
We make the following observation: the contribution to (\ref{aftertrace}) of the terms not satisfying
\begin{equation}
\label{ranges} 1\le c, m \le q^{\epsilon} \ \ \text{ and } \ \ q^{2-\epsilon}\le n \le q^{2+\epsilon}
\end{equation}
is certainly less than $q^{\epsilon}$. To see this, we may assume by (\ref{v}) that  $n < q^{2+\epsilon}m^{-2}$. Now if (\ref{ranges}) is not satisfied then
\begin{equation}
\frac{4\pi\sqrt{n} }{cq} < q^{-\epsilon},
\end{equation}
so that
\begin{equation}
J_{2k-1}\Big( \frac{4\pi\sqrt{n} }{cq}  \Big) \ll q^{-10}
\end{equation}
by (\ref{decay}) provided that $k$ is large enough.
This implies the claim. For the terms that do satisfy (\ref{ranges}), we analyze the right hand side of (\ref{aftertrace}) in dyadic intervals of $n$. To this end, let $U(x)$ be a smooth function, compactly supported on $1\le x \le 2$. By (\ref{mult}) and the observation above, we see that Theorem \ref{meanval} follows from
\begin{lemma} \label{last} For any integers
\begin{equation}
1\le c, d, m \le q^{\epsilon}, \ \ \ \ \ q^{2-\epsilon}\le N \le q^{2+\epsilon},
\end{equation}
we have that
\begin{align}
\label{sum1} \sum_{\substack{n\ge 1\\ (n,q)=1}} A_f(n,1) S(nd,1;cq) W\Big( \frac{n}{N} \Big)\ll q^{2+\epsilon},
\end{align}
where
\begin{equation}
W(x) = J_{2k-1}\Big( \frac{4\pi\sqrt{xdN} }{cq}  \Big)V\Big( \frac{ xNd^3m^2}{ q^2}\Big)U(xd)
\end{equation}
is supported on $q^{-\epsilon} \le x \le q^{\epsilon}$ and satisfies
\begin{equation}
\label{wnice} W^{(\ell)}(x) \ll_{\ell} q^{\ell \epsilon }.
\end{equation}
\end{lemma}

\section{Proofs of Lemmas \ref{voronoi} and \ref{last}}

Writing $d_1=d/(d,c)$, $c_1=c/(d,c)$, and
\begin{equation}
S(nd,1;cq) = \sums_{h \bmod cq} e\Big(\frac{ndh+\overline{h}}{cq}\Big),
\end{equation}
where $*$ indicates that the sum is restricted to primitive residue classes, we have that (\ref{sum1}) is equivalent to
\begin{align}
\label{equiva} \sums_{h \bmod cq} e\Big(\frac{\overline{h}}{cq}\Big) \sum_{\substack{n\ge 1\\ (n,q)=1}} A_f(n,1) e\Big(\frac{nd_1h}{c_1q}\Big) W\Big( \frac{n}{N} \Big)\ll q^{2+\epsilon}.
\end{align}
We now concentrate our efforts on the inner $n$-sum. The proof of Lemma \ref{voronoi} can be gleaned from this analysis on taking $d_1=c_1=1$. Right at the end, the outer $h$-sum will be executed to complete the proof of Lemma \ref{last}.

\subsection{Dirichlet characters}

Grouping the left hand side of (\ref{equiva}) by the value of $(n,c_1)$, say $d_2$, and writing $c_2=c_1/d_2$, we have that it equals
\begin{align}
 \sums_{h \bmod cq} e\Big(\frac{\overline{h}}{cq}\Big)  \sum_{d_2| c_1} \sum_{\substack{n\ge 1\\ (n,c_2q)=1}} A_f(nd_2,1)  e\Big(\frac{nd_1h}{c_2q}\Big)W\Big( \frac{nd_2}{N} \Big).
\end{align}
Note that $(nd_1h,c_2q)=1$ so that the following identity holds:
\begin{equation}
e\Big(\frac{nd_1h}{c_2q}\Big) = \frac{1}{\varphi(c_2q)}\sum_{\chi \bmod c_2q} \tau(\chi) \overline{\chi}(nd_1h),
\end{equation}
where $\varphi$ is Euler's totient function.
Writing the characters above in terms of the primitive characters which induce them and using \cite[Lemma 3.1]{iwakow}, we have that
\begin{multline}
e\Big(\frac{nd_1h}{c_2q}\Big) =\frac{1}{\varphi(c_2q)} \sum_{c_3|c_2} \Bigg( \summ_{\chi \bmod c_3q} \mu(c_2/c_3) \chi(c_2/c_3) \tau(\chi) \overline{\chi}(nd_1h)\\+\summ_{\chi \bmod c_3} \mu(c_2q/c_3) \chi(c_2 q/c_3) \tau(\chi) \overline{\chi}(nd_1h)\Bigg).
\end{multline}
Thus we have shown that to prove Lemma \ref{last} it is enough to establish that
\begin{equation}
\sums_{h \bmod cq} e\Big(\frac{\overline{h}}{cq}\Big)  \summ_{\chi \bmod c_3q} \chi(\overline{hd_1}c_2/c_3) \tau(\chi)  \sum_{\substack{n\ge 1\\ (n,c_2/c_3)=1}} A_f(nd_2,1)  \overline{\chi}(n)W\Big( \frac{nd_2}{N} \Big) \ll q^{3+\epsilon},
\end{equation}
where all the new parameters $c_i$ and $d_i$ are less than $q^{\epsilon}$. We simplify this a little more. We write
\begin{equation}
h= h_1c +h_2q,
\end{equation}
where $h_1$ varies over the primitive residue classes modulo $q$ and  $h_2$ varies over the primitive residue classes modulo $c$. We also write
\begin{equation}
\chi = \chi_1 \chi_2
\end{equation}
where $\chi_1$ varies over the primitive characters modulo $q$ and $\chi_2$ varies over the primitive characters modulo $c_3$. Thus, using \cite[(3.16),(12.20)]{iwakow}, it suffices to prove that
\begin{equation}
\label{sum2} \sums_{h_1 \bmod q} e\Big(\frac{\overline{h_1c}}{q}\Big)  \summ_{\chi_1 \bmod q} \chi_1(\overline{h_1d_1}c_2) \tau(\chi_1)  \sum_{\substack{n\ge 1\\ (n,c_2/c_3)=1}} A_f(nd_2,1)  \overline{\chi_1\chi_2}(n)   W\Big( \frac{nd_2}{N} \Big) \ll q^{3+\epsilon}.
\end{equation}

\subsection{Functional equations}\label{sect:fetwist}

To analyze the innermost sum in (\ref{sum2}), we need the functional equation of $L(s,\sym f \times \overline{\chi_1\chi_2})$.
\begin{lemma}
Let $\chi_1$ and $\chi_2$ be primitive characters mod $q$ and mod $c_3$ respectively, for $c_3<q^{\epsilon}$. Let
\begin{equation}
\alpha=\chi_1\chi_2(-1).
\end{equation}
For  some complex number $\varepsilon$, depending on $q$ and $\chi_2$, and some some integers $b_1$ and $b_2$ satisfying
\begin{equation}
|\varepsilon|=1, \ \ \ 1\le b_1, b_2 \ll  q^{\epsilon},
\end{equation}
we have that
\begin{equation}
\label{fetwist} H_{\alpha}(s) L(s, \sym f\times \overline{\chi_1\chi_2}) =   \frac{ \varepsilon  i^{\frac{1+\alpha}{2}} \overline{\chi_1}(b_1) \tau(\overline{\chi_1})^3}{q^\frac{3}{2}}  (q^{\frac{3}{2}})^{1-2s} H_{\alpha}(1-s) L(1-s, \sym f\times \chi_1\chi_2),
\end{equation}
where
\begin{equation}
\label{g2} H_{\alpha}(s)= \Gamma\Big(\frac{s+1-\frac{1+\alpha}{2}}{2}\Big) \Gamma(s+k-1) 2^{-s-1} \pi^{-\frac{3}{2}s-\frac{1}{2}} b_2^s.
\end{equation}
If $c_3=1$ (so that $\chi_2=1$), then $\varepsilon=b_1=b_2=1$. The left hand side of (\ref{fetwist}) is an entire function.
\end{lemma}
\proof

The final statement of the lemma follows from the work of Shimura \cite[Theorem 1, Theorem 2 and the following remarks]{shi}. The functional equation is a result of the work of Li \cite{li} and Atkin and Li \cite{atkli}. Specifically, in \cite[Theorem 2.2]{li}, we set
\begin{equation}
F_1= f_{\overline{\chi_1\chi_2}}
\end{equation}
to be the twist of $f$ by $\overline{\chi_1\chi_2}$,
\begin{equation}
F_2=f,
\end{equation}
and we read off the functional equation of
\begin{equation}
\label{feli} L_{F_1, F_2}(s)= L(2s, \overline{\chi_1^2\chi_2^2})\sum_{n\ge 1} \frac{\lambda_f(n)^2 \overline{\chi_1 \chi_2}(n)}{n^s}.
\end{equation}
The functional equation involves the Atkin-Lehner pseudo-eigenvalues
$\lambda_p(F_1)$ for $p|c_3q$. These are given by \cite[Theorem 4.1]{atkli}, after observing that in this theorem, the calculation of $\lambda_q(F_\chi)$ holds for $q\nmid N$ by the same proof. We get that
\begin{equation}
(2\pi)^{-2s} \Gamma(s) \Gamma(s+k-1) L_{F_1, F_2}(s) 
= A(s) (2\pi)^{-2+2s} \Gamma(1-s) \Gamma(-s+k) L_{\overline{F_1}, \overline{F_2}}(1-s),
\end{equation}
where
\begin{equation}
A(s)=\varepsilon \overline{\chi_1}(b_1) \frac{\tau(\overline{\chi_1})^4}{q^2}  \Big(\frac{4b_2q^2}{\pi^2}\Big)^{1-2s}
\end{equation}
and $\varepsilon, b_1$ and $b_2$ are as in the statement of the lemma.

By \cite[(0.4)]{shi}, we have that
\begin{equation}
L(s, \sym f\times \overline{\chi_1 \chi_2}) = \frac{L_{F_1, F_2}(s)}{L(s, \overline{\chi_1\chi_2})}. 
\end{equation}
Thus (\ref{fetwist}) follows from (\ref{feli}) and the functional equation of $L(s,\overline{\chi_1\chi_2})$, which may be found in \cite[Theorem 4.15]{iwakow} for example.
\endproof

\subsection{Summation}

Let
\begin{equation}
\widetilde{W}(s)= \int_{0}^{\infty} W(x) x^{s-1} \ dx 
\end{equation}
denote the Mellin transform of $W$, which satisfies
\begin{equation}
\widetilde{W}(s) \ll_{\Re(s)} q^{\ell \epsilon } (|s|+1)^{-\ell}
\end{equation}
for any integer $\ell\ge 0$ by (\ref{wnice}) and integration by parts $\ell$ times.
For a prime $p$ and integer $r\ge 0$, let
\begin{equation}
R_{p,r}(s) = \sum_{j=r}^{\infty} \frac{A_f(p^j,1)\overline{\chi_1\chi_2}(p^j) }{p^{js}}
\end{equation}
and
\begin{equation}
R(s) = \prod_{\substack{p|\frac{c_2}{c_3}\\p\nmid d_2}} R_{p,0}(s)^{-1}  \prod_{p^{r} \parallel d_2}  R_{p,0}(s)^{-1}  R_{p,r}(s) .
\end{equation}
Note that
\begin{align}
 \prod_{p^{r} \parallel d_2}  R_{p,0}(s)^{-1}  R_{p,r}(s) &= \prod_{p^{r} \parallel d_2}  R_{p,0}(s)^{-1}  \Big(R_{p,0}(s) -  \sum_{j=0}^{r-1} \frac{A_f(p^j,1)\overline{\chi_1\chi_2}(p^j) }{p^{js}}\Big)\\
\label{hnote1} &=  \prod_{p^{r} \parallel d_2} \Big( 1- R_{p,0}(s)^{-1}  \sum_{j=0}^{r-1} \frac{A_f(p^j,1)\overline{\chi_1\chi_2}(p^j) }{p^{js}}\Big)
\end{align}
and 
\begin{equation}
\label{hnote2} R_{p,0}(s)^{-1} = 1-\frac{A_f(p,1)\overline{\chi_1\chi_2}(p)}{p^s}+ \frac{A_f(p,1)\overline{\chi_1\chi_2}(p^2)}{p^{2s}}- \frac{\overline{\chi_1\chi_2}(p^3)}{p^{3s}}
\end{equation}
by \cite[(0.2)]{shi}.

We have that
\begin{align}
\label{melsum0} \sum_{\substack{n\ge 1\\ (n,c_2/c_3)=1}} A_f(nd_2,1)  \overline{\chi_1\chi_2}(n)   W\Big( \frac{nd_2}{N} \Big) &= \frac{1}{2\pi i}  \int_{(2)}   \sum_{\substack{n\ge 1\\ (n,c_2/c_3)=1}} \frac{A_f(nd_2,1)  \overline{\chi_1\chi_2}(n)}{n^s} \Big(\frac{N}{d_2}\Big)^{s} \widetilde{W}(s) \ ds \\
\label{melsum} &= \frac{1}{2\pi i} \int_{(2)}  R(s) L(s, \sym f \times \overline{\chi_1\chi_2})  \Big(\frac{N}{d_2}\Big)^{s} \widetilde{W}(s) \ ds .
\end{align}
By (\ref{fetwist}), we have that the integral in (\ref{melsum}) equals
\begin{equation}
 \label{melsum2} \frac{ \varepsilon  i^{\frac{1+\alpha}{2}} \overline{\chi_1}(b_1) \tau(\overline{\chi_1})^3}{q^\frac{3}{2}}  \frac{1}{2\pi i}  \int_{(2)}  R(s)   L(1-s, \sym f \times \chi_1\chi_2 )  (q^{\frac{3}{2}})^{1-2s} \Big(\frac{N}{d_2}\Big)^{s} \frac{H_{\alpha}(1-s)}{H_{\alpha}(s)} \widetilde{W}(s) \ ds.
 \end{equation}
 Keeping in mind the last statement of Lemma \ref{fetwist}, moving the line of integration to $\Re(s)=-\sigma$ for any $\sigma>0$, and using (\ref{twist}), we have that (\ref{melsum2}) equals
 \begin{equation}
\label{melsum3}  \varepsilon  i^{\frac{1+\alpha}{2}} \overline{\chi_1}(b_1) \tau(\overline{\chi_1})^3 \sum_{n\ge 1} \frac{A_f(n,1) \chi_1\chi_2(n) }{n}  \frac{1}{2\pi i}  \int_{(\sigma)} \Big( \frac{nN}{q^3d_2} \Big)^{-s} R(s)  \frac{H_{\alpha}(1-s)}{H_{\alpha}(s)} \widetilde{W}(s) \ ds.
 \end{equation}
Now, expanding out $R(s)$ and using (\ref{heckemult}, \ref{gl3coeff}, \ref{hnote1},\ref{hnote2}), we have that (\ref{melsum3}) equals
\begin{equation}
\label{melsum4} \sum_{a_1,a_2} \pm \overline{\chi_1\chi_2}(a_1) \lambda_f(a_2)  \varepsilon i^{\frac{1+\alpha}{2}}  \overline{\chi_1}(b_1) \tau(\overline{\chi_1})^3 \sum_{n\ge 1} \frac{A_f(n,1) \chi_1\chi_2(n) }{n}  \mathcal{W}_\alpha\Big(\frac{n}{q^3 N^{-1} a_1^{-1} d_2} \Big)
  \end{equation}
for a sum over some integers $a_1,a_2<q^{\epsilon}$, where 
 \begin{equation}
 \mathcal{W}_\alpha(x) =  \frac{1}{2\pi i}  \int_{(\sigma)} x^{-s}  \frac{H_{\alpha}(1-s)}{H_{\alpha}(s)} \widetilde{W}(s) \ ds
 \end{equation}
 is defined for any $x,\sigma>0$ and satisfies
 \begin{equation}
 \mathcal{W}_\alpha(x)\ll_\sigma q^{\epsilon}x^{-\sigma}.
 \end{equation}
 Thus the $n$-sum in (\ref{melsum4}) is essentially supported on $n< q^{1+\epsilon}$.
 
 \subsection{Character sums}
 
 Putting (\ref{melsum4}), the evaluation of the left hand side of (\ref{melsum0}), back into (\ref{sum2}) and exchanging the order of summation, we see that it is enough to prove that
 \begin{multline}
 \sums_{h_1 \bmod q} e\Big(\frac{\overline{h_1c}}{q}\Big) \sum_{n\ge 1} \frac{A_f(n,1)\chi_2(n) }{n}  \mathcal{W}_\alpha\Big(\frac{n}{q^3 N^{-1} a_2^{-1} d_2} \Big)  \summ_{\substack{\chi_1 \bmod q \\ \chi_1\chi_2(-1)=\alpha}} \chi_1(\overline{h_1d_1a_1 b_1}c_2n) \tau(\chi_1) \tau(\overline{\chi_1})^3 \\ \ll q^{3+\epsilon}
 \end{multline}
for $\alpha=\pm 1$ and any $a_1,a_2<q^{\epsilon}$. For $(n,q)=1$, the innermost sum equals
\begin{multline}
  \frac{q}{2} \summ_{\chi_1 \bmod q}  \chi_1(-\overline{h_1d_1a_1 b_1}c_2n) \tau(\overline{\chi_1})^2(\chi_1\chi_2(-1)\alpha + 1) \\= \frac{q(q-1)}{2}\Big(S(\overline{-h_1d_1a_1 b_1}c_2n,1;q) + \alpha\chi_2(-1)S(\overline{h_1d_1a_1b_1}c_2n,1;q)\Big)-\frac{q(1+ \alpha\chi_2(-1))}{2},
\end{multline}
while it equals 0 if $q|n$.
Thus, it remains to prove that
\begin{equation}
\label{finalstep} \sums_{h_1 \bmod q} e\Big(\frac{\overline{h_1c_2}}{q}\Big) \sum_{\substack{n\ge 1\\(n,q)=1}} \frac{A_f(n,1)\chi_2(n) }{n}  \mathcal{W}_\alpha\Big(\frac{n}{q^3 N^{-1} a_2^{-1} d_2} \Big) S(\pm \overline{h_1d_1a_1}b_1c_2n,1;q) \ll q^{1+\epsilon}.
\end{equation}
Exchanging summation again, it is easy to show that
\begin{equation}
\sums_{h_1 \bmod q} e\Big(\frac{\overline{h_1c}}{q}\Big)S(\pm \overline{h_1d_1a_1}b_1c_2n,1;q) \ll q.
\end{equation}
Now (\ref{finalstep}) follows from the immediate bound
\begin{equation}
\sum_{\substack{n\ge 1\\(n,q)=1}} \Big| \frac{A_f(n,1)\chi_2(n) }{n}  \mathcal{W}_\alpha\Big(\frac{n}{q^3 N^{-1} a_2^{-1} d_2} \Big)\Big| \ll q^{\epsilon}.
\end{equation}

\bibliographystyle{amsplain}

\bibliography{l4level6}

\end{document}